\documentclass[12pt,a4paper]{article}
\usepackage{amsmath,amsfonts,amssymb,amsthm}
\usepackage{fullpage}

\title{Dirac operators on all Podle\'s quantum spheres\\[10pt]}

\author{
Ludwik D\c{a}browski, \ Francesco D'Andrea,   \\[10pt]
Scuola Internazionale Superiore di Studi Avanzati,\\
      Via Beirut 2-4, I-34014 Trieste, Italy \\[25pt]
Giovanni Landi,  \ Elmar Wagner  \\[10pt]
Dipartimento di Matematica e Informatica,
      Universit\`a di Trieste,\\
      Via Valerio 12/1, I-34127 Trieste\\[15pt]
}

\date{22 January 2007}

\long\def\symbolfootnote[#1]#2{\begingroup%
\def\thefootnote{\fnsymbol{footnote}}\footnote[#1]{#2}\endgroup}

\newtheorem{thm}{Theorem}[section]
\newtheorem{prop}[thm]{Proposition}
\newtheorem{lem}[thm]{Lemma}
\newtheorem{cor}[thm]{Corollary}
\newtheorem{df}[thm]{Definition}
\theoremstyle{definition}

\newtheorem{rem}[thm]{Remark}

\newcommand{\beq}{\begin{equation}}
\newcommand{\eeq}{\end{equation}}

\numberwithin{equation}{section}

\newcommand{\ma}[1]{\left(\!\begin{array}{cc}#1 \end{array}\!\right)}

\newcommand{\A}{\mathcal{A}}        
\newcommand{\B}{\mathcal{B}}        

\newcommand{\C}{\mathbb{C}}         

\newcommand{\co}[2]{#1_{(#2)}}      
\newcommand{\cop}{\Delta}           
\newcommand{\dn}{{\mathord{\downarrow}}} 
\renewcommand{\H}{\mathcal{H}}      

\newcommand{\W}{\mathcal{W}}
\newcommand{\V}{\mathcal{V}}
\newcommand{\M}{\mathcal{M}}

\newcommand{\inner}[1]{\left<#1\right>}

\newcommand{\half}{{\mathchoice{\oh}{\oh}{\shalf}{\shalf}}} 
\DeclareMathOperator{\id}{id}       
\newcommand{\I}{\mathcal{I}}        
\newcommand{\ket}[1]{|#1\rangle}    
\newcommand{\lt}{\triangleright}    
\newcommand{\N}{\mathbb{N}}         
\newcommand{\nn}{\nonumber}         
\newcommand{\oh}{{\tfrac{1}{2}}}    
\newcommand{\ooh}{{\tfrac{3}{2}}}   
\newcommand{\ox}{\otimes}           
\newcommand{\R}{\mathbb{R}}         
\newcommand{\set}[1]{\{\,#1\,\}}     
\newcommand{\sg}{\sigma}            
\newcommand{\shalf}{{\scriptstyle\frac{1}{2}}} 
\newcommand{\sq}{\unskip\nobreak\kern5pt\nobreak\vrule
     height4pt width4pt depth0pt}   
\newcommand{\ssesq}{{\scriptstyle\frac{3}{2}}} 
\newcommand{\sesq}{{\mathchoice{\ooh}{\ooh}{\ssesq}{\ssesq}}} 
\DeclareMathOperator{\Tr}{Tr}       
\newcommand{\U}{\mathcal{U}}        
\newcommand{\up}{{\mathord{\uparrow}}} 
\newcommand{\Z}{\mathbb{Z}}         
\def\<#1,#2>{\langle#1,#2\rangle}   

\newcommand{\SU}{\A(\mathrm{SU}_q(2))}  
\newcommand{\podl}{\A(\mathrm{S}^2_{qt})}  
\newcommand{\su}{\U_q(\mathrm{su}(2))}  
\newcommand{\lin}{{\rm span}}  

\newcommand{\hsp}{\hspace{-1pt}}
\newcommand{\hs}{\hspace{1pt}}

\newcommand{\E}{f}
\newcommand{\F}{e}
\newcommand{\K}{k}
\newcommand{\de}{\mathrm{d}}

\newcommand{\HH}{\mathcal{H}}

\newcommand{\op}{\mathrm{OP}^{-\infty}}
\newcommand{\opz}{\mathrm{OP}^0}

\begin{document}

\maketitle

\thispagestyle{empty}

\begin{abstract}\noindent
We construct spectral triples on all Podle\'s quantum spheres ${S}^2_{qt}$.  
These noncommutative geometries are equivariant for a left action of $\su$ and are regular, even and of metric dimension~$2$. They are all isospectral to the undeformed round geometry of the sphere $S^2$. There is also an equivariant real 
structure for which both the commutant property and 
the first order condition for the Dirac operators 
are valid up to infinitesimals of arbitrary order.  
\end{abstract}

\symbolfootnote[0]
{

}

\newpage
\section{Introduction \label{sec:intro}} 

We report on further explorations of the land where noncommutative 
geometry meets quantum groups and describe even spectral triples 
$(\podl,\H,D,\gamma)$ on all the quantum 2-spheres $\podl$ introduced 
by Podle\'s \cite{Podles}. Here the real deformation parameter 
is taken to be in the interval $0<q<1$. The additional parameter $t\in[0,1]$ 
labels, among other things, the classical points, i.e., the 
1 dimensional representations of the algebra $\podl$.  

The spectral triples of the present paper generalize the $0$-dimensional 
one of \cite{DS} for the standard 2-sphere and the isospectral 
$2$-dimensional one of \cite{DLPS} for the equatorial 2-sphere. 
In particular, we have spectral triples 
of metric dimension $2$ on all Podle\'s quantum 2-spheres with the eigenvalues of the Dirac operator depending 
linearly  on an angular momentum label (and with the `correct' multiplicities). 
The Dirac operators have the same spectrum as the 
Dirac operator of the `round' metric on the usual $2$-sphere 
and we are constructing isospectral deformations.
Moreover, all these triples are regular with simple and 
discrete dimension spectrum $\Sigma=\{1,2\}$, as expected from the commutative case. 
On the standard sphere we have additional families of spectral triples with eigenvalues  of $|D|$ growing not faster that $q^{-l}$ for large $l$.  

There is a crucial equivariance of the representation of the algebra $\podl$ 
on the Hilbert space of spinors $\H$ and of the Dirac operator $D$
under a left action of the quantum enveloping algebra $\su$.
We have also an equivariant real structure of the kind `up to compact 
operators' of \cite{DLPS}, that is, such that both the `commutant property' 
and the `first order condition' for a real spectral triple \cite{ConnesReal} 
are obeyed only up to infinitesimals of arbitrary order. This phenomenon has 
also been observed in \cite{DLSSV} for the manifold of $\SU$.  Its occurrence 
for all Podle\'s quantum 2-spheres is another indication that it may be 
a characteristic feature of noncommutative manifolds of quantum groups and associated quantum (homogeneous) spaces, 
at least when asking for some degree of equivariance. 
For the standard sphere, with the eigenvalues 
of $|D|$ behaving like $| d_l | \sim q^{-l}$, the conditions on the real structure  need not be modified \cite{DS}. 

It should be mentioned that we do not study the additional `discrete' series of quantum $2$-spheres in \cite{Podles} which have algebras of  matrices as coordinate algebras.

\section{Preliminary definitions}  \label{sec:def}

In this section, we set out the basic notions concerning equivariant spectral triples. 
We start by recalling the notion of a finite summable spectral triple \cite{CBook}.

\begin{df}                                                               \label{st}
A \emph{spectral triple} $(\A,\HH,D)$ is given by a complex unital $*$-algebra $\A$, 
a faithful $*$-re\-pre\-sen\-ta\-tion $\pi$ of $\A$ by bounded operators
on a (separable) Hilbert space $\HH$, and a self-adjoint
operator $D$ (the Dirac operator) such that
\begin{itemize}
\item[(i)] $(D+\mathrm{i})^{-1}$ is a compact operator, 
\item[(ii)] $[D,\pi(a)]$ is bounded for all $a\in\A$.
\end{itemize}
With $0<\mu<\infty$, the spectral triple is said to be $\mu^+$-summable, 
if the operator $(D^2+1)^{-1/2}$ is in the Dixmier ideal 
${\mathcal{L}}^{{\mu+}}(\H)$. 
\end{df}
\noindent
We shall also call $\mu$ the \emph{metric dimension} of the triple.

A spectral triple is called \emph{even} if there exists a $\Z_2$-grading operator
$\gamma$ on $\H$, $\gamma=\gamma^*$, $\gamma^2=1$, such that the Dirac operator 
is odd and the algebra is even, i.e.,
\begin{equation}\label{gammaD}
\gamma D=-D\gamma ,  \qquad 
\pi(a)\gamma=\gamma \pi(a) , \quad  \ a\in\A.
\end{equation}

Recall \cite{ConnesReal}  that a real structure on a spectral triple $(\A,\H,D)$ 
should be given by an antiunitary operator $J$ on $\H$
fulfilling the conditions $J^2=\pm 1$, $JD=\pm DJ$ and
\begin{equation}                                     \label{eq:aJb=0}
[\pi(a), J\pi(b)J^{-1}]=0, \quad\
[[D,\pi(a)],J\pi(b)J^{-1}]=0, \quad\  a,b\in\A.
\end{equation}
It was suggested in \cite{DLPS} that one should  modify these conditions in order to
obtain a nontrivial spin geometry on the coordinate algebra 
of quantum groups and of associated 
quantum (homogeneous) spaces.  Following the lines of \cite{DLPS}, 
we impose the weaker assumption 
that \eqref{eq:aJb=0} holds only modulo infinitesimals of arbitrary high order.
Here, a compact operator $A$ is regarded as an infinitesimal of arbitrary high
order if its
singular values $s_k(A)$ satisfy $\lim_{k\to\infty}k^p s_k(A)=0$ for all
$p>0$. Therefore, throughout this paper, 
we shall use the following working definition of a real structure. 

\begin{df}                                           \label{df:real}
A \emph{real structure} $J$ on  
a spectral triple $(\A,\H,D)$ is given by an antiunitary 
operator $J$ on $\H$ such that
\begin{align}                                     \label{eq:aJb=I}
& J^2=\pm 1, \qquad JD=\pm DJ, \nn \\
& [\pi(a),J\pi(b)J^{-1}]\in\I,\qquad [[D,\pi(a)],J\pi(b)J^{-1}]\in\I , \quad  \ a,b\in\A ,
\end{align} where $\I$ is an operator ideal of infinitesimals of arbitrary high order. 
We shall name the datum $(\A,\HH,D,J)$ a \emph{real spectral triple 
(up to infinitesimals)}. 

If $(\A,\H,D,\gamma)$ is even and if in addition  
\begin{equation}
J\gamma=\pm \gamma J , 
\end{equation}
we call the datum $(\A,\HH,D,\gamma,J)$ a \emph{real even spectral triple 
(up to infinitesimals)}.
\end{df}
The signs above depend on the (metric) dimension. 
We are only interested in the case when the dimension is $2$. Then 
$J^2=- 1$, $JD= DJ$ and $J\gamma=- \gamma J$. For brevity, 
we shall drop the annotation ``up to infinitesimals'' in the sequel. 

We turn now to symmetries that will be implemented by an action of a Hopf $*$-algebra.  
In the classical case of a $G$-homogeneous spin$^c$ structure on a manifold $M$, 
the symmetry 
Hopf $*$-algebra is given by the universal enveloping algebra $\U(\mathfrak{g})$ of the Lie algebra 
$\mathfrak{g}$ of $G$. 
As we shall see, this approach will force us to consider unbounded 
$*$-re\-pre\-sen\-ta\-tions. 
For this, let $\V$ be a dense linear subspace of a Hilbert space $\H$ with inner product 
$\langle \cdot,\cdot\rangle$,  and let $\U$ be a $*$-algebra.  
An (unbounded) $*$-re\-pre\-sen\-ta\-tion of $\U$ on $\V$ is a homomorphism 
$\lambda: \U\rightarrow\rm{End}(\V)$
such that $\langle \lambda(h)v,w\rangle=\langle v,\lambda(h^*)w\rangle$ for all
$v,w\in \V$ and all $h\in\U$.

Next, let $\U= (\U, \Delta, S, \varepsilon)$ be a Hopf $*$-algebra and 
$\A$ be a left $\U$-module $*$-algebra, i.e., 
there is a left action $\lt$ of a  $\U$ on $\A$ satisfying 
\begin{equation}                                      \label{eq:mod-alg} h
\lt xy =(\co{h}{1}\lt x)(\co{h}{2}\lt y),\quad h \lt
1=\varepsilon(h)1,\quad
(h \lt x)^* = S(h)^* \lt x^* ,
\end{equation}
for all $h\in\U$ and $x,y\in\A$. As customary, we use the notation 
$\Delta(h) = \co{h}{1} \otimes \co{h}{2}$.

A $*$-re\-pre\-sen\-ta\-tion $\pi$ of $\A$ on $\V$ is called
$\U$-equivariant if there exists a $*$-re\-pre\-sen\-ta\-tion
$\lambda$ of $\U$
on $\V$ such that
\begin{equation*}                              
\lambda(h)\,\pi(x)\hs \xi = \pi(\co{h}{1} \lt x)\,\lambda(\co{h}{2})\hs \xi ,
\end{equation*}
for all $h \in \U$, $x \in \A$ and $\xi \in \V$. 
Given $\U$ and $\A$ as above,
the left crossed product $*$-algebra $\A\rtimes \U$ is defined as the
$*$-algebra generated by the two $*$-subalgebras $\A$ and $\U$ with crossed
commutation relations
\begin{equation*}
    hx=(\co{h}{1} \lt x)\co{h}{2},\quad h \in \U,\ x \in \A.
\end{equation*}
Thus $\U$-equivariant representations of $\A$ correspond to
$*$-re\-pre\-sen\-ta\-tions of $\A\rtimes \U$.

A linear operator $D$ defined on $\V$ is said to be equivariant if 
it commutes with $\lambda(h)$, i.e.,  
\begin{equation}
D\lambda(h)\hs \xi=\lambda(h)D\hs \xi ,
\end{equation} 
for all $h\in\U$ and $\xi\in\V$.  
On the other hand, an antiunitary operator $J$ is said to be equivariant 
if it leaves $\V$ invariant and if it is the antiunitary part in the polar 
decomposition of an antilinear (closed) operator $T$ that satisfies the condition
\begin{equation}
T\lambda(h)\hs \xi=\lambda(S(h)^*)T\hs \xi ,
\end{equation} 
for all $h\in\U$ and $\xi\in\V$, where $S$ denotes the antipode of $\U$.  

We collect all these equivariance requirements by giving the following 
definition~\cite{SitarzEqvt}.

\begin{df}
Let $\U$ be a Hopf $*$-algebra, $\A$ a left $\U$-module $*$-algebra and $\pi$ a 
$\U$-equi\-var\-iant representation of $\A$ on (a dense linear subspace of) a Hilbert 
space $\HH$. 
A (real even) spectral triple $(\A,\HH,D,\gamma,J)$ is called \emph{equivariant} if the 
operators $D$, $\gamma$ and $J$ are equivariant in the above sense. 
\end{df}

In the remainder of this section, we  recall 
a few analytic properties of spectral triples~\cite{CM}. 
With the operator $D$ of a spectral triple $(\A,\HH,D)$, one defines an unbounded 
derivation $\delta$ on $\B(\HH)$ by 
$$
\delta(a)=[|D|,a] , \qquad  a\in\B(\HH) .
$$
\begin{df}
The spectral triple $(\A,\HH,D)$ is said to be  \emph{regular} 
if the algebra generated by $\A$ and the commutators $[D, \A]$ belongs to 
$$
\opz:=\bigcap\nolimits_{j\in\N_0}\mathrm{dom}\,\delta^j .
$$
\end{df}
For a regular spectral triple, 
the algebra $\Psi^0$ generated 
by $\A$,  the commutators $[D, \A]$ and iterated applications of $\delta$
is a subalgebra of $\opz$,
$\Psi^0\subset\opz$.
If the spectral triple is of dimension $\mu$, the ``zeta-type'' functions 
$$
\zeta_a(z):=\Tr_{\HH}\bigl(a|D|^{-z}\bigr), \qquad 
a\in\Psi^0, 
$$ 
are defined and holomorphic for $z\in\C$ with $\mathrm{Re}\,z>\mu$.
Here we are assuming, for simplicity, that $D$ is invertible. Analogous formul{\ae} are easily defined for the general case: one works with $ (1+D^2)^{-z/2}$ instead of $|D|^{-z}$ for a not invertible $D$.

\begin{df}
A spectral triple has \emph{dimension spectrum} $\Sigma\subset\C$, with $\Sigma$ 
a countable set, if all $\zeta_a(z)$, $a \in\Psi^0$,
extend to meromorphic functions on $\C$ with poles in $\Sigma$ as unique singularities.
\end{df}

For later use, we also recall the definition of ``smoothing operators'' $\op$, i.e.,
$$
\op:=\{T\in\opz\,:\,|D|^nT\in\opz\ \;\mbox{for\ all}\ n\in\N_0\}.
$$
The class  $\op$ is a two-sided $*$-ideal in the $*$-algebra $\opz$ 
and $\delta$-invariant.
When computing the dimension spectrum, we can take the quotient 
by smoothing operators since they merely contribute holomorphic terms.

If the metric dimension of $D$ is finite, smoothing operators are infinitesimals of
arbitrary high order. But the converse is in general not true, since infinitesimals
of arbitrary high order are a two-sided $*$-ideal in $\B(\HH)$ while $\op$ is not.

Also, for a finite metric dimension, rapid decay matrices 
(in a basis of eigenvectors for $D$ 
with eigenvalues in increasing order) are smoothing operators.

\section{The Podle\'s spheres and their symmetries}       \label{sec:alg-defns}

Recall that $q\in(0,1)$ and $t\in[0,1]$. Let 
\begin{equation*}
[n] = [n]_q := \frac{q^n - q^{-n}}{q -q^{-1}},\qquad n\in\R.
\end{equation*}
The coordinate $*$-algebra $\podl$ of the Podle\'s
spheres
\cite{Podles} is generated by  elements $x_{-1}$, $x_0$ and $x_{1}$ with
relations
\begin{subequations}
\begin{eqnarray}
&x_{-1}(x_0-t)=q^2(x_0-t)x_{-1},& \\
& x_{1}(x_0-t)=q^{-2}(x_0-t)x_{1},&\\
&-[2]x_{-1}x_1 + (q^2 x_0 + t )(x_0-t ) = [2]^2 (1- t ),&\\
&-[2]x_1x_{-1} + (q^{-2} x_0 + t )(x_0-t ) = [2]^2 (1- t ),
\end{eqnarray}
\end{subequations}
and with involution $x_{-1}^*=-q^{-1}x_1$ and $x_0^*=x_0$.

The standard quantum sphere corresponds to the value $t=1$,
the equatorial one to $t=0$. For $t = 1$, the quantum sphere has one classical point. 
That is,   
there exists exactly one
1-dimensional irreducible representation of the algebra, namely, 
\beq\label{cpointt1}
x_0 = 1 , \qquad x_1=x_{-1}=0 .
\eeq
For $t \neq1$, the classical points make up an $\mathrm{S}^1$ given by the following
1-dimensional irreducible representations of the algebra,
\beq\label{cpointtnot1}
x_0 = t , \qquad x_1= \sqrt{(1+q^2)(1-t)} \ \lambda , \qquad     
x_{-1}=- q^{-1}\sqrt{(1+q^2)(1-t)} \ \bar{\lambda} , 
\eeq
where $\lambda$ is a complex number of modulus 1.

The generators that we use are related to those of \cite{Podles} as follows. 
For $t\neq 0$,  let $c=t^{-1}-t$. Setting
$$
A=(1+q^2)^{-1}(1-t^{-1}x_0), \qquad B=q(1+q^2)^{-1/2}t^{-1}x_{-1},
$$
we recover the
generators from \cite[Equation (7a)]{Podles} satisfying the relations
$$
AB=q^{-2} BA,\ \  AB^\ast=q^2 B^\ast A,\ \ B^\ast B= A-A^2+c,\ \
BB^\ast =q^2 A-q^4 A^2+c.
$$
For $t=0$, we set
$$
A=-(1+q^2)^{-1}x_0, \qquad B=q(1+q^2)^{-1/2}x_{-1},
$$
and 
obtain the generators from \cite[Equation (7b)]{Podles} with relations
$$
AB=q^{-2} BA,\ \  AB^\ast = q^2B^\ast A,\ \ B^\ast B= -A^2+1,\ \
BB^\ast=-q^4 A^2+1.
$$

\noindent
The symmetry that we shall use in the following is given by the Hopf $*$-algebra $\su$.
This Hopf $*$-algebra is generated by elements $\E$, $\F$, $\K$, $\K^{-1}$
with defining relations (cf.~\cite{KlimykS})
\begin{align*}
\K\K^{-1}=\K^{-1}\K=1, \quad
\K^2 - \K^{-2} = (q - q^{-1})(\E\F - \F\E), \quad \K\E = q\E \K,  \ \, \K\F = q^{-1}\F \K,
\end{align*}
coproduct
$$
\cop \K = \K \ox \K,  \quad
\cop \E = \E \ox \K + \K^{-1} \ox \E,  \quad
\cop \F = \F \ox \K + \K^{-1} \ox \F,
$$
counit 
$$
\varepsilon(\K-1)=\varepsilon(\E)=\varepsilon(\F)=0,
$$
antipode
$$
S(\K) = \K^{-1}, \quad S(\E)= -q \E, \quad S(\F)= -q^{-1} \F,
$$ 
and involution 
$$
\K^* = \K, \quad \E^* = \F.
$$

The irreducible finite dimensional $*$-re\-pre\-sen\-ta\-tions
of $\su$ are labeled by non-negative half-integers 
$l\in \half\N_0$ called spin.
The spin $l$ representation, say $\sg_l$, acts on a 
Hilbert space $\V^l$ which is $(2l+1)$-dimensional; the action of the generators on an orthonormal basis $\set{v^l_m : m = -l, -l+1, \cdots, l}$ 
is given by the formul{\ae},
\begin{align}
\sg_l(\K)\,v^l_m &= q^m \,v^l_m,
\nn \\
\sg_l(\E)\,v^l_m &= [l-m]^{1/2}\,[l+m+1]^{1/2} \,v^l_{m+1},
                                                 \label{eq:uqsu2-repns} \\
\sg_l(\F)\, v^l_m &= [l-m+1]^{1/2}\,[l+m]^{1/2} \, v^l_{m-1}.
\nn
\end{align}

Finally, there is  a left $\su$-action on $\podl$ 
which turns the latter into a $\su$-module
$*$-algebra. On generators, it is given by 
\begin{align*}
\K&\lt x_{-1} = q^{-1}x_{-1}, & \E&\lt x_{-1} = [2]^{1/2} x_{0}, &
\F&\lt x_{-1} = 0,\\
\K&\lt x_0 = x_0, & \E&\lt x_0 = [2]^{1/2} x_{1}, &
\F&\lt x_0 = [2]^{1/2} x_{-1},\\
\K&\lt x_1 = q x_1, & \E&\lt x_1 = 0, &
\F&\lt x_1 = [2]^{1/2} x_{0}.
\end{align*}
Note that the generators $x_m, m=-1,0,+1$, of $\podl$ transform like the basis vectors $v^1_m$
of the representation $\sg_1$.

\section{Equivariant representations} \label{sec:eqvt-repns}        

As mentioned above, looking for $\su$-equivariant representation of the algebra $\podl$, we can equivalently look for representations of the algebra $\podl\rtimes\su$; these were constructed in \cite{{SchmuedgenWCPAPod}}.

\subsection{Integrable representations of {\boldmath$\podl\rtimes\su$} } \label{sec:eqvt-repns-podl}

A $*$-re\-pre\-sen\-ta\-tion of $\podl\rtimes\su$ is called {\it integrable} if
its restriction to $\su$ is a direct sum of spin $l$ representations
$\sg_l$.

The integrable representations of
$\podl\rtimes\su$ were completely
classified in \cite{SchmuedgenWCPAPod}. It turned out that each such  
integrable representation 
is a direct sum of irreducible ones. The following proposition restates
\cite[Theorem 4.1]{SchmuedgenWCPAPod} in the case of irreducible
representations. In doing so, we replace the parameter $r$ there 
by $t^{-1}-t$.

\begin{prop}                                         \label{irreps}
Let $N\in\half\Z$. Any irreducible integrable representation of the 
$*$-algebra  $\podl\rtimes\su$ is unitarily equivalent to a representation 
$\pi_N$ described as follows. 
The representation space is $\M_N:=\bigoplus_{l=|N|,|N|+1,\cdots}\V^l$, 
where $\V^l$ is a spin $l$-re\-pre\-sen\-ta\-tion space.
The Hilbert space is the closure of $\M_N$. The generators 
of $\su$ act on each $\V^{l}$ by $\sigma_l$ as in (\ref{eq:uqsu2-repns}). 
The action of the generators $x_1$, $x_0$, $x_{-1}$ of
$\podl$ is determined on an orthonormal basis 
$\{ v^{l}_{m,N}\,:\, l=|N|,|N|+1, \cdots ,\ m=-l, -l+1, \cdots ,l\}$ 
of $\M_N$ by
\begin{equation}                                      
\label{pi}
\pi_N(x_i) v^l_{m,N} =
\alpha^-_i(l,m;N)\hspace{1pt} v^{l-1}_{m+i,N} + 
\alpha^0_i(l,m;N)\hspace{1pt} v^{l}_{m+i,N} + 
\alpha^+_i(l,m;N)\hspace{1pt} v^{l+1}_{m+i,N},
\end{equation}                                      
where the coefficients are explicitly given by
\begin{subequations}\label{eq:coeff}
\begin{align}
\alpha^+_1(l,m;N) &=  
   q^{-l+m} [l\!+\!m\!+\!1]^{1/2} [l\!+\!m\!+\!2]^{1/2}
   [2l\!+\!1]^{-1/2}[2l\!+\!2]^{-1/2}\alpha_N(l\!+\!1),     \nonumber\\
\alpha^0_1(l,m;N) &=
   -q^{m+2} [l\!-\!m]^{1/2} [l\!+\!m\!+\!1]^{1/2} [2]^{1/2} [2l]^{-1}
   \beta_N (l) ,                                                 \label{x1}\\
\alpha^-_1(l,m;N) &=  
   -q^{l+m+1} [l\!-\!m\!-\!1]^{1/2} [l\!-\!m]^{1/2}[2l\!-\!1]^{-1/2}[2l]^{-1/2}
   \alpha_N (l),  \nonumber 
\end{align}
\begin{align}
\hspace{12pt}\alpha^+_0(l,m;N) &=
   q^m [l\!-\!m\!+\!1]^{1/2} [l\!+\!m\!+\!1]^{1/2}[2]^{1/2}
   [2l\!+\!1]^{-1/2}[2l\!+\!2]^{-1/2}\alpha_N(l\!+\!1),     \nonumber\\
\alpha^0_0(l,m;N) &=           
   [2l]^{-1} \big([l\!-\!m\!+\!1][l\!+\!m] - q^{2} [l\!-\!m][l\!+\!m\!+\!1]\big)
   \beta_N (l) ,                                               \label{x0}  \\
\alpha^-_0(l,m;N) &=           
   q^m [l\!-\!m]^{1/2} [l\!+\!m]^{1/2} [2]^{1/2}[2l\!-\!1]^{-1/2} [2l]^{-1/2}
   \alpha_N (l),   \nonumber
\end{align}
\begin{align}
\alpha^+_{-1}(l,m;N) &=
   q^{l+m} [l\!-\!m\!+\!1]^{1/2}[l\!-\!m\!+\!2]^{1/2}[2l\!+\!1]^{-1/2} 
   [2l\!+\!2]^{-1/2}\alpha_N (l\!+\!1),           \nonumber\\               
\alpha^0_{-1}(l,m;N) &=
   q^m [l\!-\!m+1]^{1/2} [l\!+\!m]^{1/2} [2]^{1/2}
   [2l]^{-1} \beta_N (l) ,                                       \label{x-1}\\
 \alpha^-_{-1}(l,m;N) &=
   -q^{-l+m-1} [l\!+\!m\!-\!1]^{1/2} [l\!+\!m]^{1/2}[2l\!-\!1]^{-1/2}[2l]^{-1/2} 
   \alpha_N (l)  \nonumber
\end{align}
\end{subequations}
(with the convention that $\alpha^0_{i}(0,0;0), \alpha^{-}_{i}(0,0;0),
\alpha^-_{i}(\half,\pm\half;\pm\half)$ are zero).
The real numbers $\beta_N(l)$ and $\alpha_N(l)$ are given by
\begin{subequations}
\begin{align}                                           
\label{alpha+-}
\alpha_N(l) &=
\frac{[2]^{1/2}[l+N]^{1/2}[l-N]^{1/2}[2l]^{1/2} }{[2l+1]^{1/2}[l]}
\Big(1-t+q^{-2N}[2l]^{-2}[l]^{2}(t-1+q^{2N})^2\Big)^{1/2}, \\
\beta_N(l) &=
\frac{
\epsilon \, \hs [\, 2|N|\,
](q^{-1}\hsp + \hsp q \hsp - \hsp q^{\epsilon}\hs t ) + t (q\hsp -\hsp
q^{-1})\{[|N|]\hspace{1pt}[|N|\!+\!1] - [l][l\!+\!1]\}
}{q[2l+2]}
\label{beta+-},
\end{align}
\end{subequations}
with $\epsilon={\rm sign}(N)$.
\end{prop}
\begin{rem}
In \cite[Proposition 6.4]{SchmuedgenWCPAPod}, it was shown that
the representation space $\M_N$
is isomorphic to a projective left $\podl$-module of rank 1. These
projective left $\podl$-modules
can be considered as line bundles over the quantum sphere
$\mathrm{S}^2_{qt}$ with winding numbers $2N$ \cite{BM,HM,MS}.
In particular, $\M_0\cong \podl$ is the trivial line bundle.
\end{rem}
\begin{rem}
The representations $\pi_{\pm\half}$ for the standard quantum sphere 
$\A(\mathrm{S}^2_{q1})$ appeared first in \cite{DS}.
\end{rem}

\subsection{The spin representation} \label{sec:spin-repns-podl}

Classically, the spinor bundle on the 2-sphere is given by the direct sum
of two line bundles with winding numbers $-1$ and 1.
In correspondence with
the classical picture, we take 
\[
\W:=\M_{-\half}\oplus\M_{\half}
\]
as (polynomial sections) of the quantum spinor bundle.
The Hilbert space completion $\H$ of $\W$ will be the
{\em Hilbert space of spinors} with the representation 
$\pi:=\pi_{-\half}\oplus\pi_{\half}$ of the algebra $\podl$ as the {\em spinor representation}.

For simplicity of notation,  we shall use the label $+$ 
instead of $N = \frac{1}{2}$, and $-$ instead of $N = -\frac{1}{2}$.
With these conventions, $\W:=\M_{-}\oplus\M_{+}$, and $\pi_\pm$  denotes the 
representation $\pi_{\pm\half}$. 
Moreover we shall identify elements of  
$\su$ with the corresponding operators on $\W$ coming from Proposition~\ref{irreps}. By that proposition we have the decomposition
$$
\W =\bigoplus_{l=\frac{1}{2},\frac{3}{2},\cdots}\W^l, \quad \
\W^l=\bigoplus_{m=-l,\cdots,l} \W^l_m,
$$
where
$\W^l_m=\lin\{v^l_{m,-}, v^l_{m,+} \}$. Throughout this paper, we keep the
basis
$\{v^l_{m,-}, v^l_{m,+} \}$ of $\W^l_m$ fixed.  
Clearly, an orthonormal basis of $\H$ is provided by the vectors
\begin{equation}                                      
\label{basis}
\{v^l_{m,-},\  v^l_{m,+}\,:\,l=\tfrac{1}{2},\tfrac{3}{2},\cdots,\ m=-l,\cdots,l\} .
\end{equation}                                      

For an arbitrary but fixed vector
$w^l_m\in \W^l_m$,
say $w^l_m=\alpha^l_{m,-}v^l_{m,-}+ \alpha^l_{m,+} v^l_{m,+}$ with components
$\alpha^l_{m,-}, \alpha^l_{m,+}\in\C$, and for any $i,j\in\Z$, the expression $w^{l+j}_{m+i}$ denotes the
vector $w^{l+j}_{m+i}=\alpha^l_{m,-}\hs v^{l+j}_{m+i,-}+ \alpha^l_{m,+}\hs
v^{l+j}_{m+i,+}\in \W^{l+j}_{m+i}$. It is understood
that $w^{l+j}_{m+i}=0$ whenever $l+j<0$ or $|m+i|>l+j$.

A substantial part of our results will be based on the fact that the operators
$\pi(x_{-1})$, $\pi(x_0)$, $\pi(x_1)$ can be ``approximated''
by operators acting  diagonally on $\M_{-}\oplus \M_{+}$.
This is the content of the next lemma.
For this purpose, we define operators $z_{-1}$, $z_0$ and $z_1$ 
on $\H$ by their action on $w^l_m\in \W^l_m$ as follows:
\begin{equation}                                      \label{eq:appr}
z_iw^l_{m}=\alpha^-_i(l,m;0)\hspace{1pt} w^{l-1}_{m+i}
+\alpha^0_i(l,m;0)\hspace{1pt} w^{l}_{m+i}
+\alpha^+_i(l,m;0)\hspace{1pt} w^{l+1}_{m+i}\hs,\quad i=-1,0,1 .
\end{equation}
Here, the coefficients are the ones defined in Equations~\eqref{eq:coeff} 
unless $|m+i|>l+\nu$,\,\ $\nu=0,\pm 1$; in this case, we set $\alpha^\nu_i(l,m;0)=0$.
Formally, the operators $z_i$'s are given by the same formul{\ae} as
the $\pi_0(x_i)$'s but $\W$ and $\M_0$ decompose into different spin representation spaces.  
As a consequence, $z_{-1}$, $z_0$ and $z_1$ do not satisfy the commutation relations  of the generators of 
$\podl$. Explicitly, their commutation relations read
\begin{align*}
z_{-1}(z_0-t)&=q^2(z_0-t)z_{-1}, \\
z_{1}(z_0-t)&=q^{-2}(z_0-t)z_{1}, \\
-[2]z_{-1}z_1 + (q^2 z_0 + t )(z_0-t )&= [2]^2 (1- t )-cP_-\hs , \\
-[2]z_1z_{-1} + (q^{-2} z_0 + t )(z_0-t )&= [2]^2 (1- t )-cP_+\hs ,
\end{align*}
where $P_\pm$ denotes the projection onto the $2$-dimensional
subspace $\,\lin\{v^{1/2}_{\pm 1/2,-}, v^{1/2}_{\pm 1/2,+} \}$ and 
$\,c:=q[2]\alpha_0(\half)^2=q[2]\big(1-t+[\frac{1}{2}]^2t^2\big)$.

Before stating the lemma, we need some more notation. 
For any bounded operator $L$ on the Hilbert space $\H$,\, $\I(L)$ denotes 
the two-sided ideal of $\B(\H)$ generated by $L$. 
Also, throughout this paper, $L_q$ stands for 
the compact operator on $\H$ defined by
\beq
L_q w^l_{m}=q^l w^l_{m}, \quad w^l_m\in \W^l_m .
\eeq

\begin{lem}                                                   \label{appr}
There exist bounded operators $A_i$ and $B_i$,
$i=-1,0,1$, such that
\begin{align*}
&\pi(x_i)-z_i=A_i L_q^2=L_q^2 B_i  \qquad  \mbox{if}\,\ t\neq 1,\quad \\
& \pi(x_i)-z_i=A_i L_q=L_q B_i  \qquad  \mbox{if}\,\ t=1.
\end{align*}
In particular, $\pi(x_i)-z_i\in\I(L_q^2)$ for $t\neq 1$ and
$\pi(x_i)-z_i\in\I(L_q)$ for $t=1$ 
\end{lem}

\begin{proof}
Our aim is to show that $L_q^{-k}(\pi(x_i)-z_i)$ and
$(\pi(x_i)-z_i)L_q^{-k}$ are bounded, with $k=2$ for $t\neq 1$
and $k=1$ for $t=1$. Notice that $\pi(x_i)-z_i$ is a sum of three
independent weighted shift operators with weights
$\alpha^\nu_i(l,m;\pm\frac{1}{2})-\alpha^\nu_i(l,m;0)$, where $\nu\in\{+,-,0\}$.
Hence it suffices to prove that
$q^{-kl}|\alpha^{\nu}_i(l,m;\pm\frac{1}{2})-\alpha^{\nu}_i(l,m;0)|\leq C^{\nu}_i$
for some constants $C^{\nu}_i\in\R^+$.
Using the inequalities 
\begin{equation}
\label{ineq}
[n]\leq (q^{-1}-q)^{-1} q^{-n}\;, ~\forall n\geq 0\, , \qquad
 [n]^{-1} \leq q^{n-1}\;, ~\forall n\geq 1\, ,
\end{equation}
one verifies that the  coefficients in front of
$\alpha_{N}(l)$ and $\beta_{N}(l)$ in Equations \eqref{x1}-\eqref{x-1} are uniformly
bounded. Thus, we only have to prove that the sequences
$q^{-kl}|\beta_N(l)-\beta_{0}(l)|$ and $q^{-kl}|\alpha_N(l)-\alpha_{0}(l)|$ are bounded.
We prove this for $N\in \half\Z$, even though we need it only for $N=\pm 1/2$.
The boundedness of the first sequence follows from the identity
$$
\bigl|\beta_N(l)-\beta_0(l)\bigr|=
\frac{\big|\epsilon[2|N|](q^{-1}+q-q^{\epsilon}t)+t(q-q^{-1})[|N|][|N|+1]\big|}{q[2l+2]}\;,
$$
using $(q[2l+2])^{-1}\leq q^{2l}$.
Concerning the second sequence, if $t=1$ one verifies directly that
for any $N$ the sequence $q^{-l}|\alpha_N(l)|$ is bounded;
thus $q^{-l}|\alpha_N(l)-\alpha_{0}(l)|$ is bounded.
If $t\in[0,1)$, set 
$$
u_l:=\frac{[l+N][l-N][2l]}{q[2l+1][l]^2}\,
\frac{1-t+q^{-2N}[2l]^{-2}[l]^{2}(t-1+q^{2N})^2}{1-t}-1 \;.
$$
for $l \geq 1$. From $\,[l+N][l-N][2l]\leq q[2l+1][l]^2\,$ we obtain
\begin{equation*}
u_l\leq
\frac{q^{-2N}[2l]^{-2}[l]^2(t-1+q^{2N})^2}{1-t}=
C_N(1-q^2)^2[2l]^{-2}[l]^2\leq C_Nq^{2l} \;,
\end{equation*}
with $C_N:=\frac{\{q^N-q^{-N}(1-t)\}^2}{(1-q^2)^2(1-t)}$.
On the other hand from $[l+N][l-N]=[l]^2-[N]^2$ and $[2l]-q[2l+1]=-q^{2l+1}$ we get
$$
u_l\geq\frac{[2l][l+N][l-N]}{q[2l+1][l]^2}-1=
\frac{-[2l][N]^2-q^{2l+1}[l]^2}{q[2l+1][l]^2}\geq
-q^{2l-2}\frac{[N]^2}{1-q^2}-q^{4l-1} \;.
$$
Hence there exists a constant $C'_N\in\R^+$ such that $|u_l|\leq C'_Nq^{2l}$.
Moreover from the inequality $|\sqrt{1+u}-1|\leq |u|$ which holds 
for any $u\geq -1$, it follows that
$$
\left|\frac{\alpha_N(l)}{\sqrt{(1+q^2)(1-t)}}-1\right|=
\bigl|\sqrt{1+u_l}-1 \bigr|\leq |u_l|\leq C'_Nq^{2l}\;.
$$
Finally, using the triangle inequality we conclude that
$\frac{|\alpha_N(l)-\alpha_0(l)|}{\sqrt{(1+q^2)(1-t)}}\leq\,(C_0+C_N)q^{2l}$.
\end{proof}

\section{Equivariant Dirac operators}  \label{sec:Dirac}

In this section, we provide a class of equivariant self-adjoint operators $D$ 
on the Hilbert space of spinors $\H$ such that 
$(\podl,\H,D)$ fulfills all conditions of a  spectral triple.

Let  $C_q:=q\K^2+q^{-1}\K^{-2}+(q-q^{-1})\F\E$ be  the Casimir operator of
$\su$. It acts on the space of spinors $\W$ by $C_q\hs
v^l_{m,\pm}=(q^{2l+1}+q^{-2l-1})v^l_{m,\pm}$,
and its closure on  $\H$ is self-adjoint.
Suppose that $D$ is a self-adjoint operator on $\H$ commuting strongly
with the closure of $C_q$. Then the finite-dimensional subspaces $\W^l$
reduce $D$, and $\W$ is invariant under $D$. Assume  next  that $D$ is
equivariant, i.e.,
$XD=DX$ for all $X\in\su$. 

From $\K D=D\K$, it follows that $D$ leaves the
2-dimensional
subspaces $\W^l_m$ invariant. Let $D^l_m$ be a self-adjoint operator on
$\C^2$ such
that the restriction of $D$ to $\W^l_m$, in our fixed
basis for the latter, is given by $D^l_m$.
Then $D\E=\E D$ implies that $D^l_{m+1}=D^l_m$. Hence there exists a
self-adjoint operator on $\C^2$, $D_l$ say, such that $D^l_m=D_l$ for all
$m=-l,\cdots, l$. Diagonalizing $D_l$, we can write
\begin{equation}                                       \label{defD}
D_l= U_l^*
\begin{pmatrix}
d_l^\up & 0\\
0 & d_l^\dn
\end{pmatrix} U_l,
\end{equation}
where $d_l^\up$ and $d_l^\dn$ are the eigenvalues of $D_l$ and $U_l$ 
is a unitary operator on $\C^2$. Without loss of generality,
we may assume that $d_l^\dn\leq d_l^\up$. Definition \ref{st} (i) 
is fulfilled  if and only if
\begin{equation}                                       \label{compact}
\lim_{l\rightarrow\infty}|d_l^\up| =\lim_{l\rightarrow\infty}|d_l^\dn|
=\infty.
\end{equation}

Next we give sufficient conditions for the boundedness of
the   commutators $[D,\pi(a)]$ with $a\in\podl$. 
By the Leibniz rule for the commutator,  $[D,xy]=[D,x]y+x[D,y]$, it suffices to
consider only
commutators with the generators $x_{-1}$, $x_0$, $x_1$.
As the action
of the operators $z_{-1}$, $z_0$, $z_1$ from Lemma~\ref{appr}
is comparatively simple,
we write
$[D,\pi(x_i)]\hsp =\hsp [D,z_i]\hsp +\hsp [D,\pi(x_i)\hs {-}\hs z_i]$ 
and start by analyzing $[D,z_i]$.
Let $w^l_m\in \W^l_m$. Then
\begin{align}
[D,z_i]w^l_m=&\,\alpha^-_i(l,m;0)\Big( U_{l-1}^*
\begin{pmatrix}
d_{l-1}^\up & 0\\
0 & d_{l-1}^\dn
\end{pmatrix} U_{l-1} -
U_l^*\begin{pmatrix}
d_l^\up & 0\\
0 & d_l^\dn
\end{pmatrix} U_l\Big)w^{l-1}_{m+i} \nonumber\\
& + \alpha^+_i(l,m;0)
\Big( U_{l+1}^*\begin{pmatrix}
d_{l+1}^\up & 0\\
0 & d_{l+1}^\dn
\end{pmatrix} U_{l+1} -
U_l^*\begin{pmatrix}
d_l^\up & 0\\
0 & d_l^\dn
\end{pmatrix} U_l\Big)w^{l+1}_{m+i}.
\end{align}
We need to treat the cases $t\neq 1$ and $t=1$ separately.

Firstly, let $t\neq 1$. Using \eqref{ineq}
and observing that 
$\lim_{l\rightarrow\infty}\alpha_{\pm\frac{1}{2}}(l)=(1+q^2)^{1/2}(1-t )^{1/2}$
we deduce that the sequence $\{ \alpha^\pm_i(l,m;0) \}$ 
is uniformly bounded and does not converge to zero. 
Hence $[D,z_i]$ is bounded if and only if there exist $C\in\R$ such that
\begin{equation}                                        \label{D}
\big|\big|
U_{l+1}^*\begin{pmatrix}
d_{l+1}^\up & 0\\
0 & d_{l+1}^\dn
\end{pmatrix} U_{l+1} -
U_l^*\begin{pmatrix}
d_l^\up & 0\\
0 & d_l^\dn
\end{pmatrix} U_l
\big|\big| < C\quad \mbox{for all} \ l=\half,\sesq,\cdots \, , 
\end{equation}
or, equivalently,
\begin{equation}                                        \label{UDU-D}
\big|\big|
 U_l U_{l+1}^*\begin{pmatrix}
d_{l+1}^\up & 0\\
0 & d_{l+1}^\dn
\end{pmatrix} U_{l+1}U_l^* -
\begin{pmatrix}
d_l^\up & 0\\
0 & d_l^\dn
\end{pmatrix}
\big|\big| < C\quad \mbox{for all} \ l=\half,\sesq,\cdots \,.
\end{equation}
Denoting by $t_l$ the absolute value of the (1,1) entry of the  matrix
$U_l U_{l+1}^*$
the last inequality implies that there is 
a fixed constant $C^\prime\in\R$ satisfying
\begin{align*}
|d_{l+1}^\up-d_{l}^\up-t_l^2(d_{l+1}^\up-d_{l+1}^\dn)|
=|d_{l+1}^\dn-d_{l}^\up+(1-t_l^2)(d_{l+1}^\up-d_{l+1}^\dn)| & <C^\prime , \\  
|d_{l+1}^\dn-d_{l}^\dn+t_l^2(d_{l+1}^\up-d_{l+1}^\dn)|
=|d_{l+1}^\up-d_{l}^\dn-(1-t_l^2)(d_{l+1}^\up-d_{l+1}^\dn)| & <C^\prime , \\
t_l\sqrt{1-t_l^2}\hs |d_{l+1}^\up-d_{l+1}^\dn| & <C^\prime , 
\end{align*}
for all $l=\half,\sesq,\cdots$.
Using $t_l\in[0,1]$ and $ d_{l+1}^\up\geq d_{l+1}^\dn$, one
deduces that
$$
\sup\{\,|d_{l+1}^\dn-d_{l}^\dn|,\,|d_{l+1}^\up-d_{l}^\up|\,:\,
l=\half,\sesq,\cdots\,\}<\infty.
$$
Consequently, there exist $C\in\R$ such that
\begin{equation}                                \label{beschr}
|d_{l}^\dn|<C\hs l,\quad
|d_{l}^\up|<C\hs l.
\end{equation}
From these inequalities, it follows that both 
$DL_q^2$ and $L_q^2D$ are bounded and so is
$[D,\pi(x_i)-z_i]$ by Lemma~\ref{appr}. 
Hence Equation  \eqref{D} implies that both commutators
$[D,z_i]$ and $[D,\pi(x_i)-z_i]$ are bounded and so is their sum 
$[D,\pi(x_i)]$.

If $t=1$, a direct computation shows that the sequence
$\{q^{-l}\alpha^\pm_i(l,m;0) \}$
is uniformly bounded and does not converge to zero. As a consequence, $[D,z_i]$ is bounded if and
only if there are $C\in\R$ such that
\begin{equation}                                             \label{D0q}
\big|\big|
U_{l+1}^*\begin{pmatrix}
d_{l+1}^\up & 0\\
0 & d_{l+1}^\dn
\end{pmatrix} U_{l+1} -
U_l^*\begin{pmatrix}
d_l^\up & 0\\
0 & d_l^\dn
\end{pmatrix} U_l
\big|\big| < C q^{-l}\quad \mbox{for all} \ l=\half,\sesq,\cdots \,.
\end{equation}
On the other hand, a sufficient condition for the commutator $[D,\pi(x_i)-z_i]$ to be  bounded is that $DL_q=L_qD$ is bounded (since 
by Lemma~\ref{appr}, $\pi(x_i)-z_i=A_i L_q=L_q B_i$ with bounded
operators $A_i$ and $B_i$). 
Now the operator $DL_q=L_qD$ is
bounded if and only if there exist $C\in\C$ such that
\begin{equation}                                           \label{dl0}
|d_{l}^\dn|<Cq^{-l},\quad |d_{l}^\up|<C\hs q^{-l}, \quad
l=\half,\sesq,\cdots \,.
\end{equation}
Since these inequalities imply \eqref{D0q}, they are   sufficient for the commutator 
$[D,\pi(x_i)]=[D,z_i]+[D,\pi(x_i)-z_i]$ to be bounded.

\bigskip
Summarizing, we have established the following proposition.
\begin{prop}\label{Dirac} 
Let
$D$ be an equivariant self-adjoint operator on $\H$ such that $D$ and the
closure of $C_q$ strongly commute. Then
$\W^l_m=\lin\{v^l_{m,-}, v^l_{m,+} \}$ is an invariant subspace and the
restriction of $D$ to $\W^l_m$ is given by \eqref{defD}. 
The triple $(\podl,\H,D)$ defines an $\su$
equivariant spectral triple if, for $t\neq 1$, conditions \eqref{compact} and \eqref{D} are satisfied;  and for $t=1$,  conditions \eqref{compact} and \eqref{dl0} are satisfied. 
\end{prop}

\section{The real structure} \label{sec:real}

We require that the real structure  $J$ be equivariant.
This means that $J$ is the antiunitary part of a closed antilinear operator
$T$  satisfying $TX=S(X)^*T$ for all $X\in\su$
(cf.~Section~\ref{sec:def}). Now, any antiunitary operator $J$ 
which leaves $\W$ invariant and fulfills
\begin{equation}                                         \label{JXJ}
      JXJ^{-1}=\K S(X)^*\K^{-1},     \qquad  
X\in\su,
\end{equation}
is equivariant. To see this, one can take for $T$ the closure of the operator $J\K$ 
since $J\K X=S(X)^*J\K$ for all $X\in\su$.

Next, consider the antiunitary operator $J_0$ on $\H$ defined by
\begin{equation}
 J_0\hs v^l_{m,\pm} = (-1)^{m+1/2}\hs v^l_{-m,\pm}, \quad
l=\half,\sesq,\cdots,\quad
m=-l,\cdots,l.
\end{equation}
Clearly, $J_0^2=-1$.
The antiunitary operator $J_0$ will play a crucial role in discussing the
general form of an equivariant real structure.
We summarize some properties of
$J_0$ in the following lemmata.
\begin{lem}                              \label{J0lem}
The antiunitary operator $J_0$ is equivariant. 
\end{lem}
\begin{proof}
First note that $J_0$ leaves the spinor bundle $\W$ invariant. 
The lemma is proved by showing that $J_0$ satisfies \eqref{JXJ}.
Since $X\mapsto \K S(X)^*\K^{-1}$ is
an antilinear homomorphism, it suffices to verify \eqref{JXJ}  for the
generators of $\su$ which can easily be done by straightforward
calculations. \end{proof}
\begin{lem}
\label{ziJzj} With the operators $z_{-1}$, $z_0$, $z_1$ defined in 
Equation  \eqref{eq:appr}, 
\begin{equation}                         \label{ziJ0zj}
      [z_i, J_0 z_j J_0^{-1}]=0, \quad i,j=-1,0,1.
\end{equation}
\end{lem}
\begin{proof}
In the notation of Equation  \eqref{eq:appr}, we have
\begin{equation}
\label{eq:apprJ0}
J_0z_iJ_0^{-1}w^l_{m}=(-1)^{i}\big(\alpha^-_{i}(l,-m;0)\hspace{1pt}
w^{l-1}_{m-i}
+\alpha^0_{i}(l,-m;0)\hs w^{l}_{m-i}
+\alpha^+_{i}(l,-m;0)\hs w^{l+1}_{m-i}\big) ,
\end{equation}
for $w^l_m\in \W^l_m$.
The lemma is proved by direct computations using Equations
\eqref{eq:appr} and \eqref{eq:apprJ0}.
\end{proof}

\begin{lem}                                                   \label{aJ0b} 
Let $\I$ denotes the operator ideal $\I(L_q)$ or $\I(L_q^2)$
for
$t=1$ and $t\neq 1$, respectively. Then, for all $a,b\in\podl$,
\begin{equation}                                                \label{xiJxj}
 [\pi(a),J_0\hs \pi(b)\hs J_0^{-1}]\in\I.
\end{equation}
\end{lem}

\begin{proof}
Again from the Leibniz rule of commutators, it suffices to prove
Equation  \eqref{xiJxj} for the generators $x_{-1}$, $x_0$, $x_1$.
Since
$$
[\pi(x_i),J_0\pi(x_j)J_0^{-1}]=[\pi(x_i)-z_i,J_0\pi(x_j)J_0^{-1}]+
[z_i,J_0\big(\pi(x_j)-z_j\big)J_0^{-1}]+[z_i, J_0 z_j J_0^{-1}],
$$
the assertion follows from Lemmata \ref{appr} and \ref{ziJzj}.
\end{proof}

\begin{lem}                                         \label{J0}
Let $(\podl,\H,D,J)$ be an equivariant 2-dimensional real spectral triple. 
Assume that $J$ satisfies \eqref{JXJ}. 
Then there exist unitary operators $W_l$ on
$\C^2$ such that the restrictions of $D$ and $J$ to $\W^l_m$ can be
expressed as 
\begin{align}                                       \label{DW}
D\lceil_{\W^l_m}&=D_l= W_l^*
\begin{pmatrix}
d_l^\up & 0\\
0 & d_l^\dn
\end{pmatrix} W_l,                          \\       \label{JW}
J\lceil_{\W^l_m}&=W_l^* \overline{W}_l J_0\lceil_{\W^l_m}
=W_l^*  J_0\lceil_{\W^l_m} W_l,
\end{align}
where the bar denotes complex conjugation.
\end{lem}

\begin{proof}
Given an antiunitary operator $J$ satisfying \eqref{JXJ}, 
the  operator $V=JJ_0^{-1}=-JJ_0$
is unitary and fulfills for all $X\in\su$ 
$$
VX V^*=JJ_0XJ_0^{-1}J^{-1}=S(S(X)^*)^*=X 
$$
since, for any Hopf $*$-algebra, $*\circ S\circ *\circ S=\rm{id}$. In
particular, $V$ commutes with $\K$ and with the Casimir operator $C_q$ and does
so with their closures. By an argument similar to the one at the beginning of
Section~\ref{sec:Dirac}, one shows that there exist unitary operators
$V_l$ on $\C^2$ such that the restriction of $V$ to $\W^l_m$ is given by
$V_l$. Thus 
$$
J=VJ_0 , \qquad J\lceil_{\W^l_m}=V_l J_0\lceil_{\W^l_m} .
$$
Recall that, in our standard basis, $D_l=D\lceil_{\W^l_m}$ has the form
\eqref{defD}.
As $J_0$ is antiunitary, we get
$$J_0\lceil_{\W^l_m} D_l=\overline{D}_l J_0\lceil_{\W^l_m}
$$ 
(the bar is  still complex conjugation). 
From the requirement $JDJ^{-1}=D$, it follows that $D_l=V_l \overline{D}_l V_l^{*}$
or,
equivalently,
\begin{equation*}
\begin{pmatrix}
d_l^\up & 0\\
0 & d_l^\dn
\end{pmatrix}
= U_l V_l \overline{U}_l^*
\begin{pmatrix}
d_l^\up & 0\\
0 & d_l^\dn
\end{pmatrix} (U_l V_l \overline{U}_l^*)^*.
\end{equation*}
If $d_l^\up\neq d_l^\dn$, the last equation implies that $U_l V_l \overline{U}_l^*$
is a diagonal matrix. Thus there are angle variables $\phi_l,\psi_l\in[0,2\pi)$ such that 
\begin{equation}                                           \label{Vl}
V_l=U_l^*
\begin{pmatrix}
\exp(\mathrm{i}\phi_l) & 0\\
0 & \exp(\mathrm{i}\psi_l)
\end{pmatrix}
\overline{U}_l.
\end{equation}
Inserting
$$
W_l:=
\begin{pmatrix}
\exp(-\frac{\mathrm{i}}{2}\phi_l) & 0\\
0 & \exp(-\frac{\mathrm{i}}{2}\psi_l)
\end{pmatrix}
U_l
$$
into \eqref{defD} and \eqref{Vl}, we arrive at \eqref{DW} and \eqref{JW}.
Clearly, Equation  \eqref{DW} remains valid if $d_l^\up=d_l^\dn$ and $W_l$
is an
arbitrary unitary matrix. Therefore we may assume without loss of
generality that
$D$ and $J$ are given on $\W$ by Equations \eqref{DW} and \eqref{JW},
respectively.
\end{proof}

Under the assumption that $J$ satisfies \eqref{JXJ}, 
Lemma~\ref{J0} provides necessary conditions 
for equivariant real spectral triples. 
Next, we aim at finding sufficient conditions for $(\A,\HH,D,J)$ to yield an 
equivariant real spectral triple up to infinitesimals.
The ideal $\I$ from Definition \ref{df:real} 
is generated by operators $L_q^\beta$, where $\beta$ is a
positive real constant.  
Again, the basic idea of the proof is 
to replace the generators $\pi(x_i)$ by the simpler operators $z_i$.
By the last lemma, the operator $D$ and $J$ act diagonally with respect 
to the basis formed by  $\tilde v^l_{m,\pm}:=W_l^* v^l_{m,\pm}$, 
but the operators $z_i$ do not do so (in general). 
Their action on this basis is found to be
\begin{equation*}                                    
z_i\tilde v^l_{m,\pm}=\alpha^-_i(l,m;0)\hspace{1pt} W_lW^*_{l-1}
\tilde v^{l-1}_{m+i,\pm}
+\alpha^0_i(l,m;0)\hspace{1pt} \tilde v^{l}_{m+i,\pm}
+\alpha^+_i(l,m;0)\hspace{1pt}W_lW^*_{l+1} \tilde v^{l+1}_{m+i,\pm}\hs.
\end{equation*}
In order to reduce the computations to diagonal operators, we require thus 
that 
\begin{equation}                                          \label{Wcond} 
||W_{l} W_{l+1}^*-1||< Cq^{\beta l}
\end{equation}
for some real constant $C$. Then  \eqref{Wcond} implies that, 
modulo the ideal $\I(L_q^\beta)$, 
the operators $z_i$ act on the basis vectors $\tilde v^l_{m,\pm}$ 
as a sum of diagonal shift operators. 
Next, with the operators $z_i$ 
acting on the basis vectors $\tilde v^l_{m,\pm}$as
diagonal shifts, 
we compute 
\begin{align*}
|\langle \tilde v^{l+2}_{m+i-j,+},
[[D,z_i],J z_jJ^{-1}&]\tilde v^l_{m,+}\rangle|\\
&=\alpha^+_i(l\hsp +\hsp 1,m\hsp - \hsp j;0)\hs
\alpha^+_{j}(l,-\hsp m;0) 
|( d^\up_{l+2}  -  d^\up_{l+1})-
(d^\up_{l+1}  - d^\up_l )|,
\end{align*}
where we used  $\alpha^+_i(l\hsp +\hsp 1,m\hsp - \hsp j;0)\hs
\alpha^+_{j}(l,- m;0)=\alpha^+_{j}(l\hsp +\hsp 1,- m\hsp-\hsp i;0)
\hs\alpha^+_i(l,m;0)$. 

Let $t\neq 1$. To ensure that 
$|\langle \tilde v^{l+2}_{m+i-j,+},
[[D,z_i],J z_jJ^{-1}]\tilde v^l_{m,+}\rangle|<Cq^{\beta l}$ 
for some real constant $C$, we must impose the condition 
\begin{equation}                                        \label{uplin}
|( d^\up_{l+2}  -  d^\up_{l+1})-
(d^\up_{l+1}  - d^\up_l )|<Cq^{\beta l} ,
\end{equation}
since $\alpha^+_i(l\hsp +\hsp 1,m\hsp - \hsp j;0)\hs
\alpha^+_{j}(l,-\hsp m;0)=O(1)$. A similar argument leads also to 
\begin{equation}                                      \label{dnlin}
|( d^\dn_{l+2}  -  d^\dn_{l+1})-
(d^\dn_{l+1}  - d^\dn_l )|<Cq^{\beta l}.
\end{equation}
From Equation~\eqref{eq:aJb=I}, Lemma~\ref{J0} and the foregoing,
it is immediately clear that one can always add to $D$ a
self-adjoint operator from $\I(L_q^\beta)$ having the same eigenvectors.
A sufficient condition for Equations \eqref{uplin} and \eqref{dnlin} 
to be satisfied is to 
assume that the eigenvalues $d^\up_l$ and $d^\dn_l$ depend linearly on $l$.
But this dependence is alike the one in \eqref{beschr} 
to get bounded commutators of $D$ with algebra elements.

On the other hand, for $t=1$, the argument leading to Equations \eqref{uplin} 
and \eqref{dnlin} fails since in this case $\alpha^+_i(l,m;0)=O(q^l)$. However, 
a linear dependence on the eigenvalues of $D$ clearly satisfies 
the condition \eqref{dl0} for bounded commutators in this case, too.  
Thus, by Proposition \ref{Dirac}, for all Podle\'s quantum spheres,
we have  equivariant spectral triples $(\podl,\H,D)$ 
when the eigenvalues 
$d^\up_l$ and $d^\dn_l$ of $D$ depend linearly on $l$. 

The next proposition shows that this linear dependence, together with the requirement \eqref{Wcond}, 
suffices to obtain equivariant real spectral triples (up to infinitesimals) 
for all quantum spheres $\podl$.

 \begin{prop}                                             \label{realDirac}
Suppose that
\begin{equation}                                           \label{dl}
d^\up_l=c^\up_1 l+c^\up_2, \quad d^\dn_l=c^\dn_1 l+c^\dn_2,
\end{equation}
where $c^\up_1$, $c^\up_2$, $c^\dn_1$, $c^\dn_2$ are real numbers such
that $c^\up_1\neq 0$ and $c^\dn_1\neq 0$. Let $D$ be the self-adjoint operator on
$\H$ determined by \eqref{DW}, and let the antiunitary operator $J$ be
given by
\eqref{JW}. Suppose that  $\beta\in (0,2]$
for $t\neq 1$ and $\beta\in (0,1]$ for $t=1$, and
assume that \eqref{Wcond} holds.
Then $(\podl,\H,D,J)$ defines an equivariant real
spectral
triple such that Equation  \eqref{eq:aJb=I} is satisfied with
$\I:=\I(L_q^\alpha)$, where $\alpha\in(0,\beta)$.
\end{prop}

\begin{proof}
Clearly, Equations \eqref{DW} and \eqref{dl} uniquely determine a
self-adjoint
operator  (denoted by $D$) on $\H$ since
the collection $\{ W_l^*v^l_{m,-},   W_l^*v^l_{m,+}\,:\,l=\half,\sesq,\cdots,\
m=-l,\cdots\hsp,l\}$
is a complete set of orthonormal eigenvectors.
By Lemma~\ref{J0}, it  is evident that
$JD=DJ$ and that $J^2=-1$. Lemma~\ref{J0lem} and Equation  \eqref{JW} imply that
\eqref{JXJ} also applies, so $J$ is equivariant.
Thus it remains to verify Equation 
\eqref{eq:aJb=I}.

Observe that $\I(L_q^{\alpha_1})\subset \I(L_q^{\alpha_2})$ for
$\alpha_2\geq\alpha_1>0$. From
\begin{align*}
(Jz_iJ^{-1}-J_0z_iJ_0^{-1})w^l_m                
&= (-1)^{i}\Big(\alpha^-_{i}(l,-m;0)\hs 
(W_{l-1}^* \bar W_{l-1} \bar W_{l}^* W_{l}-1) w^{l-1}_{m-i}  \\
& \qquad \qquad \qquad + \alpha^+_{i}(l,-m;0)\hs
(W_{l+1}^* \bar W_{l+1} \bar W_{l}^* W_{l}-1)
w^{l+1}_{m-i}\Big)
\end{align*}
and Equation \eqref{Wcond}, it follows  
that $Jz_iJ^{-1}-J_0z_iJ_0^{-1}\in \I(L_q^\beta)$.
As a consequence, and using Lemma~\ref{ziJzj}, one finds that 
\begin{equation}                                          \label{zJzJ}
[z_i, J z_j J^{-1}]=[z_i, J_0 z_j J_0^{-1}]+[z_i,Jz_iJ^{-1}-J_0 z_j
J_0^{-1}] \in \I(L_q^\beta).
\end{equation}
Now the proof
of the first equation of \eqref{eq:aJb=I} is completely analogous to the
proof of Lemma~\ref{aJ0b} with $J_0$ replaced by $J$.

By a repeated use of the Leibniz rule for the commutator and of the
first relation of \eqref{eq:aJb=I}, we need to
prove the second relation of \eqref{eq:aJb=I} only for the generators
$x_{-1}$, $x_0$, $x_1$. Note that the operator $DL_q^{\beta-\alpha}=L_q^{\beta-\alpha}
D$ is bounded. Hence, by Lemma~\ref{appr},  
it follows that $[D,\pi(x_i)-z_i]\in\I(L_q^\alpha)$. Now, writing
\begin{align*}
[[D,\pi(x_i)], & J \pi(x_j)J^{-1}]  \\
& =[[D,\pi(x_i)-z_i],J\pi(x_j)J^{-1}] 
+[[D,z_i],J(\pi(x_j)-z_j)J^{-1}]
+[[D,z_i],Jz_jJ^{-1}]
\end{align*}
and using again Lemma~\ref{appr}, we see that it suffices to establish that
\begin{equation}                                       \label{DzJz}
[[D,z_i],Jz_jJ^{-1}]\in\I(L_q^\alpha).
\end{equation}
Let $W$ denote the unitary operator on $\H$ given by $Ww^l_m=W_lw^l_m$ for
all
$w^l_m\in\W^l_m$. Then $J=W^* J_0 W$ and $D=W^*\widetilde D W$, where $\widetilde D$ 
is the unique self-adjoint operator on $\H$ such that 
\beq \label{diagD}
\widetilde D\lceil_{\W^l_m}=\widetilde D_l  
:=
\begin{pmatrix}
d_l^\up & 0\\
0 & d_l^\dn
\end{pmatrix}.
\eeq
In these terms,  the requirement \eqref{DzJz} is equivalent to
\begin{equation}                                       \label{DWzJ}
[[\widetilde D, Wz_iW^*],J_0 W z_jW^* J_0^{-1}]\in\I(L_q^\alpha).
\end{equation}
Evaluating $Wz_iW^*-z_i$ on vectors $w^l_m\in\W^l_m$ gives
$$
(Wz_iW^*-z_i)w^l_m
=\alpha^-_i(l,m)_0(W_{l-1}W_l^*-1) w^{l-1}_{m+i}
+\alpha^+_i(l,m)_0 (W_{l+1}W_l^*-1)w^{l+1}_{m+i}\hs.
$$
From this and \eqref{Wcond}, we conclude that
$Wz_iW^*-z_i\in\I(L_q^\beta)\subset\I(L_q^\alpha)$. Thus, Equation 
\eqref{DWzJ} is, in turn, equivalent to
\begin{equation}                                           \label{DzJ}
[[\widetilde D, z_i],J_0 z_j J_0^{-1}]\in\I(L_q^\alpha).
\end{equation}
Note now that $[[\widetilde D, z_i],J_0 z_j J_0^{-1}]$ 
can be written as a sum of five independent
weighted
shift operators  with weights $S^\nu_{i,j}(l,m)$,
$\nu=-2,\cdots,2$,
i.e.,
\begin{equation}                                          \label{shiftD} 
[[\widetilde D, z_i],J_0 z_j J_0^{-1}]w^l_m=\sum_{\nu=-2}^2 S^\nu_{i,j}(l,m)
w^{l+\nu}_{m+i-j},\quad
w^l_m\in \W^l_m.
\end{equation}
Moreover,
\begin{equation}                                       \label{DzzD}
[\widetilde D, z_i]v^l_{m,\pm} = c_1^\pm
\big(\alpha^+_i(l,m;0)\hspace{1pt} v^{l+1}_{m+i,\pm}
-\alpha^-_i(l,m;0)\hspace{1pt} v^{l-1}_{m+i,\pm}\big),
\end{equation}
where $c_1^+=c_1^\up$ and $c_1^-=c_1^\dn$.
Using \eqref{DzzD} and \eqref{eq:apprJ0}, a 
lengthy but straightforward computations shows that
$|S^\nu_{i,j,\pm}(l,m)|< C q^{2l}$ for some $C\in\R$. From this, 
we conclude that $[[\widetilde D, z_i],J_0 z_j J_0^{-1}]\in\I(L_q^2)$ which implies
\eqref{DzJ}, and the proof is complete. 
\end{proof}

\begin{cor}
Up to adding elements from the ideal $\I(L_q^\alpha)$, the operators $D$ and $J$
from Proposition \ref{realDirac} are given by
\begin{align}                                       \label{DW0}
D\lceil_{\W^l_m}&= W_0^*
\begin{pmatrix}
d_l^\up & 0\\
0 & d_l^\dn
\end{pmatrix} W_0,                         \\       \label{JW0}
J\lceil_{\W^l_m}&=W_0^*  J_0\lceil_{\W^l_m} W_0,
\end{align}
where $W_0$ is a unitary operator on $\C^2$.
\end{cor}

\begin{proof}
Let $D_1$ and $J_1$ be given by Equations \eqref{DW} and \eqref{JW}, and suppose
that the unitary operators $W_l$ satisfy \eqref{Wcond}. Then
$\{W_l\}_{l=\half,\sesq,\cdots }$ is a Cauchy sequence and
$W_0:=\lim_{l\to\infty}W_l$ is a unitary operator.
Let $D$ and $J$ denote the operators determined by \eqref{DW0} and
\eqref{JW0}.
With $\widetilde D_l$ the diagonal matrix defined in \eqref{diagD}, we have 
$$
(D_1-D)w^l_m=(W_l\widetilde D_l W_l-W_0\widetilde D_l W_0)w^l_m , 
$$ 
for
$w^l_m\in\W^l_m$.
By \eqref{Wcond}, there exist $C\in\R$ such that 
$||W_l-W_0||<C q^{\beta l}$. 
Furthermore, since $\alpha < \beta$, by \eqref{dl} it follows that  
$\lim_{l \to\infty}q^{(\beta-\alpha)l}|d^\up_l |
=\lim_{l \to\infty}q^{(\beta-\alpha)l}|d^\dn_l |=0$ and the sequences 
$\{q^{(\beta-\alpha)l}|d^\up_l | \}$ and $\{q^{(\beta-\alpha)l}|d^\dn_l | \}$ 
are uniformly bounded.
Thus, 
\begin{align*}
||W_l\widetilde D_l W_l-W_0\widetilde D_l W_0||
&=||(W_l\hsp -\hsp W_0)\widetilde D_l W_l+W_0\widetilde D_l(W_l\hsp -\hsp W_0)|| \\
& <2 ||W_l-W_0||\,||D_l|| < C^\prime q^{\alpha l} ,
\end{align*}
for some $C^\prime\in\R$, and $D_1-D\in\I(L_q^\alpha)$.
Similarly one
shows that $J_1-J\in\I(L_q^\alpha)$.
\end{proof}

\begin{rem} 
For the standard sphere, corresponding to $t=1$, the 
conditions \eqref{dl0} for bounded commutators of the Dirac 
operator with algebra elements
allow more than linear dependence for the eigenvalues of $D$. 
In \cite{DS}, the eigenvalues of $D$ were taken to be 
$q$-analogues of the spectrum of the classical Dirac operator 
of the round metric of the sphere $\mathrm{S}^2$; 
they behave like $| d_l | \sim q^{-l}$ for large $l$. 
For this family, one gets a spectral triple; this is also a particular case of our 
Proposition \ref{Dirac}. Moreover, one has a stronger result on the real structure: 
it is not up to infinitesimals but the stronger relations \eqref{eq:aJb=0} are satisfied.  
\end{rem}

\section{Equivariant real even spectral triple} \label{sec:grading}

In this section, an additional character enters the stage, the even structure. 
As shown in the next proposition, the existence of a grading operator determines  
completely the structure of the geometry $(\podl,\H,D,J,\gamma)$.

\begin{prop} \label{struct}
Let 
$(\podl,\H,D,J,\gamma)$ be an
equivariant real even spectral triple with real structure $J$ satisfying 
condition \eqref{JXJ}.
Then the operators $\gamma$,
$D$, $J$  act on the spinor bundle $\W$ by
\begin{align} \label{DJg}
& \gamma\hs v^l_{m,\pm}=\pm\hs \epsilon \hs\hs v^l_{m,\pm}\hs , \qquad 
J v^l_{m,\pm}=\zeta_l  \hs  \hs  (-1)^{m+1/2}\hs \hs
v^l_{-m,\mp}\hs , \nn \\ 
& D\hs v^l_{m,+}= d_l\hs\hs v^l_{m,-}\hs ,\qquad  D\hs v^l_{m,-}=\bar d_l\hs\hs v^l_{m,+}\hs ,
\end{align}
where $d_l\in\C$, $\epsilon\in\{-1,1\} $, and $\zeta_l\in\C$ such that
$|\zeta_l|=1$.
\end{prop}

\begin{proof}

Clearly, $\gamma\neq 1$ as otherwise the first condition in Equation~\eqref{gammaD} is not satisfied. From
$\gamma^*=\gamma$ and $\gamma^2=1$, it follows that $\gamma$ has
eigenvalues $\pm 1$. Since $\gamma$ commutes with all elements from the
crossed product algebra
$\podl\rtimes\su$ and
since the integrable representation of $\podl\rtimes\su$ on the spinor
bundle $\W$
decomposes into two inequivalent irreducible ones on
$\M_-$ and $\M_+$, we conclude that
$\gamma\hs v^l_{m,+}=\epsilon\hs v^l_{m,+}$ and
$\gamma\hs v^l_{m,-}=-\epsilon\hs v^l_{m,-}$, where $\epsilon\in\{-1,1\} $.

Recall from Lemma~\ref{J0} that $J\lceil_{\W^l_m}=V_l J_0\lceil_{\W^l_m}$
with a unitary operator  $V_l$ on $\C^2$. In addition, the condition $J\gamma=-\gamma J$
implies that $J$ maps $\M_\pm$ into $\M_\mp$. Thus 
$$
J v^l_{m,\pm}=\zeta_{l,\pm}  \hs J_0\hs v^l_{m,\mp} ,
$$ 
with
complex numbers $\zeta_{l,-} $ and $\zeta_{l,+}$ such that
$|\zeta_{l,-} |=|\zeta_{l,+} |=1$. From $J^2=-1$, we obtain 
$\zeta_{l,-} \bar \zeta_{l,+}=1$, so $\zeta_{l,-} =\zeta_{l,+} =:\zeta_{l}$.

Similarly, the condition $D \gamma = -\gamma D$ implies that $D$  maps $\M_\pm$ into
$\M_\mp$. 
Hence, by \eqref{defD},
$D$ has the form described in \eqref{DJg}.
\end{proof}

We  combine Propositions \ref{realDirac} and \ref{struct} to present
equivariant real even spectral triples. 

\begin{prop}                                             \label{teo}
Let  $\gamma$ and $D$ be  self-adjoint operators and $J$ be an 
antiunitary operator on $\HH$ given by 
\begin{equation}                                      \label{thm}
\gamma\hs v^l_{m,\pm}=\pm\hs  v^l_{m,\pm}\hs ,\quad
D\hs v^l_{m,\pm}= (c_1 l + c_2) v^l_{m,\mp}\hs ,
\quad J v^l_{m,\pm}= (-1)^{m+1/2}\hs \hs v^l_{-m,\mp}\hs ,
\end{equation}
where $c_1,c_2\in\R$ with $c_1\neq 0$. Then
$(\podl,\H,D,J,\gamma)$ is an
equivariant real even spectral triple such that the 
conditions \eqref{eq:aJb=I} are satisfied with
$\I=\I(L_q^\alpha)$, where $\alpha\in(0,2)$ for $t\neq 1$ and
$\alpha\in(0,1)$ for $t=1$. 
\end{prop}

\begin{proof} For $l=\half,\sesq,\cdots$,
set $d_l^\up=c_1 l + c_2$ and $d_l^\dn=-(c_1 l + c_2)$. Define
\begin{equation}\label{w0}
W_0:=
\frac{1}{\sqrt{2}}
\begin{pmatrix}
1 & 1\\
-\mathrm{i}& \mathrm{i}
\end{pmatrix}, \quad W_l:=W_0 .
\end{equation}
The restrictions of $D$ and $J$ to $\W^l_m$ are then given as in \eqref{DW} and \eqref{JW}, 
respectively. With the choice \eqref{w0}, the inequality \eqref{Wcond} holds 
trivially and we can suppose that $\beta=2$
for $t\neq 1$ and $\beta=1$ for $t=1$.
By Proposition \ref{realDirac}, $(\podl,\H,D,J)$ 
yields an equivariant real spectral
triple such that the conditions \eqref{eq:aJb=I} are satisfied with
$\I=\I(L_q^\alpha)$. It is obvious that $\gamma$ fulfills all the
requirements of a grading operator. 
\end{proof}

\begin{rem}
Recall from Subsection~\ref{sec:spin-repns-podl} that $\W$ can be
considered as a
deformation of the classical spinor bundle.
The classical spectral triple on
the commutative 2-sphere with its round metric is real and even, and the
corresponding
Dirac operator  has eigenvalues $\pm(l+\half)$,
$l=\half,\sesq,\cdots$,
each  with multiplicity $2l+1$.
Therefore we can regard
the equivariant real even spectral triple
$(\podl,\H,D,J,\gamma)$ from Theorem \ref{thm}
with $c_1=1$  and $c_2=\half$ as an isospectral deformation of the
classical spin geometry.
\end{rem}

\begin{rem}
By perturbing both $D$ and $J$ by  infinitesimals belonging to the ideal $\I(L_q^\alpha)$, one
produces more examples of equivariant real even spectral triples.
However, those from Theorem \ref{thm}
are distinguished by being obtained from the isospectral deformation
via rescaling linearly the eigenvalues of $D$.
\end{rem}

Finally, we prove the non-triviality of our noncommutative geometry.

\begin{prop}\label{cor:one}
Let $F\in\B(\HH)$ be the sign of the Dirac operator; it is given
by 
$$Fv^l_{m,\pm}=v^l_{m,\mp}.
$$ 
The datum $(\A(S^2_{qt}),\HH,F)$ is a $1^+$-summable
non-trivial Fredholm module.
\end{prop}

\begin{proof}
Since the chiral spin representations $\pi_\pm$ coincide modulo smoothing operators,
the commutator $[F,x]$ is a smoothing operator for all $x\in\A(S^2_{qt})$, thus
$[F,x]\in\mathcal{L}^1(\HH)$. This shows that $(\A(S^2_{qt}),\HH,F)$
is a $1^+$-summable Fredholm module.

As representative of the corresponding periodic cyclic cohomology class
$\mathrm{ch}^F$ we can take the cochain having only one component $\mathrm{ch}^F_0$
in degree $0$, given by,
\begin{equation}\label{eq:chzero}
\mathrm{ch}^F_0(a):=\tfrac{1}{2}\Tr(\gamma F[F,a])\;.
\end{equation}

The non-triviality of our Fredholm module is proved by pairing $\mathrm{ch}^F$ with the $K$-theory class 
of the projection
$$
p=\frac{1}{(2-t)(1+q^2)}\ma{1+q^2+x_0-tq^2 & -(1+q^2)^{1/2}x_1 \\ q(1+q^2)^{1/2}x_{-1} & 
1+q^2-q^2x_0-t}, 
$$
describing line bundles over the spheres ${S}^2_{qt}$ \cite{BM}.
The pairing is
\begin{align}\label{nontriv}
\inner{\mathrm{ch}^F,[p]} 
&:=\tfrac{1}{2}\Tr_{\HH\otimes\C^2}(\gamma F[F,p]) 
=\frac{1-q^2}{2(2-t)(1+q^2)}\Tr_{\HH}(\gamma F[F,x_0]) \nn \\
&=q^{-2}(1-q^2)^2\sum_{l,m}\frac{[l-m+1][l+m]}{[2l][2l+2]}\;\;.
\end{align}
For the last equality we have used the explicit formul{\ae} \eqref{x0} and \eqref{beta+-} for the coefficients $\alpha^0_0(l,m;\pm)$  
in the representation of $x_0$.
The last series in \eqref{nontriv} was shown in~\cite[Section~5.3]{DD} to be a continuous function in $q$ for $q\in[0,1)$.
Since it is integer valued in the interior of this interval (being the index of a Fredholm operator), 
it is constant by continuity and can be 
computed at $q=0$. The result is $\inner{\mathrm{ch}^F,[p]}=1$.
\end{proof}

\section{Analytic properties} \label{sec:analytic}

We describe now some further analytic properties of the spectral triple
given  in Proposition \ref{teo} (forgetting its real even structure)
and specify, for simplicity, the constants $c_1=1$ and $c_2=1/2$. 
With this choice,
$D$ is invertible.
For other allowed values, the subsequent results remain valid with 
minor changes in the proofs and by working 
with $ (1+D^2)^{-z/2}$ instead of 
$|D|^{-z}$ for a not invertible $D$. 
We use the polar decomposition $D = F |D|$, with $F$  the sign of $D$.

Recall that the spin representation $\pi = \pi_-\oplus \pi_+$ for each  
of the generators of the algebra $\podl$ has three terms appearing on the right hand side of (\ref{pi}).
Accordingly,  we shall write 
$\pi (x_i ) = x_i =: x_i^- +x_i^0 + x_i^+ $, $i= -1, 0, 1$.
The operators $x_i^{\nu}$ are weighted shifts, 
mapping each $W^l_m$ into $W^{l+\nu}_{m+i}$,\,\ $\nu=-1,0,1$,  
which are easily seen to be bounded.

\begin{prop}
The spectral triple $(\podl,\H,D)$ 
is $2^+$-summable and regular.
\end{prop}
\begin{proof}
The $2^+$-summability follows from the linear growth of the spectrum of $D$
(with the appropriate multiplicities).
To show that the spectral triple is regular, it is enough 
to prove that the elements $x_i^\nu$ and $[D,x_i^\nu]$, $\nu=0,\pm 1$, 
are in the smooth domain of  the derivation 
$\delta(\,\cdot\,):=[|D|,(\,\cdot\,)]$.
This is clear for $x_i^\nu$ since 
$$[|D|, x_i^\nu ] v^l_{m,\pm} = 
\big( (l+ \nu +\half ) \alpha^\nu_i(l,m;\pm) - (l +\half )\alpha^\nu_i(l,m;\pm) \big) 
v^{l+\nu}_{m+i,\pm},$$
hence $\delta(x_i^\nu) = \nu x_i^\nu$. 
As for $[D,x_i^\nu]$, we use the fact that 
$\delta^k([D, x_i^\nu ] )=[D, \delta^k(x_i^\nu )]= \nu ^k [D, x_i^\nu]$.
The boundedness of the latter follows from the formul{\ae} in 
Proposition \ref{irreps} giving, for instance,
$$[D, x_i^+ ] v^l_{m,\pm} = \Big( 
(l+ \half) \big(\alpha^+_i(l,m;\pm) - \alpha^+_i(l,m;\mp) \big) 
+ \alpha^+_i(l,m;\pm) \Big)\,
v^{l+1}_{m+i,\mp}.$$
The last summand is clearly bounded while the boundedness of the other 
follows from the fact that $\alpha^+_i(l,m;\pm) - \alpha^+_i(l,m;\mp)$ 
is at least of order $q^l$ 
by Equations \eqref{eq:coeff} and \eqref{alpha+-}.
Similar arguments work for $[D, x_i^- ]$.
\end{proof}

To compute the dimension spectrum, we introduce another representation of 
the algebra $\podl$ which is obtained from simpler operators.
Let $\hat{\HH}$ be a Hilbert space with orthonormal basis
$\ket{l,m}_{\pm}$, where $l\in\frac{1}{2}\Z$ and $l+m\in\N_0$.
Consider the bounded operators $\alpha,\beta\in\B(\hat{\HH})$ defined by
\begin{equation*}
\alpha\ket{l,m}_\pm
= \sqrt{1-q^{2\smash[t]{(l+m+1)}}}\ket{l+\tfrac{1}{2},m+\tfrac{1}{2}}_\pm,
\qquad \beta\ket{l,m}_\pm= q^{l+m}\ket{l+\tfrac{1}{2},m-\tfrac{1}{2}}_\pm.
\end{equation*}
These operators
satisfy the commutation relations of $\A(SU_q(2))$, i.e.,
\begin{equation*}
\beta\alpha=q\alpha\beta,\quad
\beta^*\alpha=q\alpha\beta^*,\quad
[\beta,\beta^*]=0,\quad
\alpha^*\alpha+q^2\beta^*\beta=1,\quad
\alpha\alpha^*+\beta\beta^*=1.
\end{equation*}
The embedding of the Podle\'s spheres into $\A(SU_q(2))$ (see e.g.~\cite{Podles}) 
leads to a $*$-re\-pre\-sen\-ta\-tion
$\varphi:\A(S^2_{qt})\to\B(\hat{\HH})$ given by
\begin{align*}
\varphi(x_1) &=\sqrt{1+q^2}\,
\Big\{\sqrt{1-t}\,\bigl(\alpha^2-q(\beta^*)^2\bigr)-t\beta^*\alpha\Big\}, \\
\varphi(x_0) &=(1+q^2)\Big\{\sqrt{1-t}\,(\alpha\beta+\beta^*\alpha^*)
-t\beta\beta^*\Big\}+t.
\end{align*}
Later on, we shall need the explicit expression 
\begin{align}
(1+q^2)^{-1}\varphi(x_0-t)\ket{l,m}_\pm = \;  &
\sqrt{1-t}\,q^{l+m-1}\sqrt{1-q^{2(l+m)}}\ket{l-1,m}_\pm 
\nonumber 
-t\,q^{2(l+m)}\ket{l,m}_\pm \\   &                                                  \label{expl}
+\sqrt{1-t}\,q^{l+m}\sqrt{1-q^{2(l+m+1)}}\ket{l+1,m}_\pm.
\end{align}
Our Hilbert space of spinors $\HH$ with the basis \eqref{basis}
is identified with a subspace
of $\hat{\HH}$, that is, we consider the inclusion
$$
Q: \HH\to\hat{\HH},\quad Q\hs v^l_{m,\pm} = \ket{l,m}_\pm ,
$$
for $l\in\N_0+\frac{1}{2}$ and $|m|\leq l$.
Let $P\hs\hs:\hs\hs \hat{\HH}\to\HH$ be the adjoint map of $Q$, i.e., 
$$
P\hs \ket{l,m}_\pm =v^l_{m,\pm} \ \ \mbox{for}\ 
l\in\N_0+\mbox{$\frac{1}{2}$}\ 
\mbox{and}\  |m|\leq l, \quad P\hs \ket{l,m}_\pm =0 \ \ \mbox{otherwise}.
$$
Our Dirac operator $D$ on $\HH$ is the ``restriction'' 
of the self-adjoint operator $D'$ on $\hat{\HH}$ determined by 
$$
D'\ket{l,m}_\pm=(l+\tfrac{1}{2})\ket{l,m}_\mp ,
$$
in the sense that 
$DP=PD'$, $QD=D'Q$. The same holds for $|D'|$ and for $F':=D'|D'|^{-1}$.
The subspace $Q\HH$ is not invariant for the representation $\varphi$. 
However, we can sandwich 
$\varphi$ between $Q$ and $P$ thus obtaining a $*$-linear map
\begin{equation}\label{tildephi}
\tilde{\varphi}:\A(S^2_{qt})\to\B(\HH), \qquad 
\tilde{\varphi}(a) =P\hs\hs \varphi(a)\hs\hs Q,  
\end{equation}
that has the following approximation property. 

\begin{lem}\label{lem:ap}
The operator $a-\tilde{\varphi}(a)$ is a smoothing operator on $\HH$ for all $a\in\A(S^2_{qt})$.
\end{lem}

\begin{proof}
Observe that if 
$T:\hat{\HH}\to\HH$ is a matrix of rapid decay (in our fixed bases), then so are 
$a T \varphi(b)$, where $a,b\in \A(S^2_{qt})$, and $TQ:\HH\mapsto \HH$. 
Using arguments similar to the ones in the proof of Lemma \ref{appr}, 
for the generators $x_{i}$, $i=-1,0,1$, one verifies by direct computations that 
$x_i P-P\varphi(x_i)$ yields a rapid decay matrix. Then the lemma follows from 
the identity
\begin{equation*}
abP-P\varphi(ab)=a\big(bP-P\varphi(b)\big)+\big(aP-P\varphi(a)\big)\varphi(b) 
\end{equation*}
by applying $PQ=\id_\HH$ and by the above observations. 
\end{proof}

The next proposition is the main result of this section.

\begin{prop}
The dimension spectrum is $\Sigma=\{1,2\}$.
\end{prop}
\begin{proof}
Let $\Psi^0$ be the algebra generated by
$\podl$, by $[D, a]$ for all $a\in\podl$ and by iterated applications
of the derivation $\delta$ (cf.~Section \ref{sec:def}). 
Let $\mathcal{C}$ be the $*$-algebra (of bounded operators on $\hat{\HH}$)
generated by $\alpha,\beta,\alpha^*,\beta^*$ and $F$. 
By Lemma~\ref{lem:ap}, $\A(S^2_{qt})\subset P\,\mathcal{C} Q+\op$. 
Note that 
\begin{align*}
[F',\alpha]&=0, & [D',\alpha]&=\tfrac{1}{2}\,\alpha \,F', &
[|D'|,\alpha]&=\tfrac{1}{2}\,\alpha, \\
[F',\beta]&=0, &
[D',\beta ]&=\tfrac{1}{2}\,\beta  \,F', &
[|D'|,\alpha]&=\tfrac{1}{2}\,\beta.
\end{align*}
Thus $P\,\mathcal{C}Q$ is invariant under application of $\delta$ and  $[D,(\cdot)]$ and hence 
$\Psi^0\subset P\,\mathcal{C} Q+\op$.

We shall compute the singularities of zeta functions associated to 
the  monomials $S:=P\alpha^n\beta^j(\beta^*)^kQ\hs F$ and 
$T:=P\alpha^n\beta^j(\beta^*)^kQ$, where  
$n\in\Z$ and $j,k\in\N_0$ and we employ the notation
$\alpha^n:=(\alpha^*)^{|n|}$ for $n<0$.
From the commutation relations of $\alpha$ and $\beta$, it is clear that these monomials
span $P\mathcal{C}Q$.

Firstly, note that the $\zeta$ function associated with a bounded off-diagonal operator is
identically zero in the half-plane $\mathrm{Re}\,z>2$
and so is its holomorphic extension to the entire complex plane.
This is the case for the monomials $S$ due to the presence of $F$. The other monomials 
$T$
shift the index $l$ by $(n+j-k)/2$ and the index $m$ by $(n-j+k)/2$ and 
therefore are also off-diagonal operators unless these shifts are zero, 
which happens when $n=0$ and $j=k$.
Hence only monomials $T=P\beta^k(\beta^*)^kQ=P(\beta\beta^*)^kQ$ contribute
to the dimension spectrum.

For $k=0$, $T=\id$, and the
corresponding zeta function is
$$
\zeta_{\id}(z)=\sum_{l+\frac{1}{2}\in\N}\sum_{l+m=0}^{2l}\sum_{r=\pm}
(l+\tfrac{1}{2})^{-z}=4\zeta(z-1),
$$
where $\zeta(z)$ is the Riemann zeta function, meromorphic in $\C$ with
a simple pole at $1$ and with residue $1$. 
Since $\id\in \Psi^0$, this shows that $2\in \Sigma$.

When $k> 0$, $T= P(\beta\beta^*)^kQ$. So 
$T v^l_{m,\pm} =q^{2k(l+m)}v^l_{m,\pm}$
and the associated zeta function is
\begin{align*}
\zeta_T(z) &=\sum_{l+\tfrac{1}{2}\in\N}2(l+\tfrac{1}{2})^{-z}
\sum_{l+m=0}^{2l}(q^{2k})^{l+m} 
=\sum_{l+\tfrac{1}{2}\in\N}2(l+\tfrac{1}{2})^{-z}\frac{1-q^{2k(2l+1)}}{1-q^{2k}} \\ &
=\tfrac{2}{1-q^{2k}}\zeta(z)+\textit{holomorphic function}.
\end{align*}
Therefore $\Sigma$ may contain, besides 2,  at most the additional point 1. 
We still have to check that $1\in \Sigma$ 
since the algebra $P\mathcal{C}Q$ is strictly larger than $\Psi^0$.
For this, we take $a\in\A(S^2_{qt})$, where
$(1+q^2)^2 a=(x_0-t)^2=\tilde{\varphi}\bigl((x_0-t)^2\bigr)\,+$ {\it smoothing terms}.
Then, using \eqref{expl}, we get (modulo holomorphic functions)
\begin{align*}
\zeta_{a}(z) &\sim 2\sum_{l+\frac{1}{2}\in\N}(l+\tfrac{1}{2})^{-z}
\sum_{k=0}^{2l}\Big\{t^2q^{4k}+(1-t)\big\{(1+q^{-2})q^{2k}-(q^2+q^{-2})q^{4k}\big\}\Big\} \\ &
\sim 2\frac{1+(1-t)^2}{1-q^4}\,\zeta(z),
\end{align*}
and
$\,\mathrm{Res}_{z=1}\zeta_{a}(z)=2\frac{1+(1-t)^2}{1-q^4}\neq 0$ 
for all  $t\in[0,1]$. This shows that $\Sigma=\{1,2\}$.
\end{proof}

Let $\A(S^1)$ denote the polynomial $*$-algebra in one variable $\lambda$, 
with $\lambda\bar\lambda=1$.
For $t\neq 1$, we have $*$-algebra morphisms
$\sigma_t:\A(S^2_{qt})\to \A(S^1)$ given by the `classical points' \eqref{cpointtnot1},
$$
\sigma_t(x_0)=t,\qquad\sigma_t(x_1)=\sqrt{(1+q^2)(1-t)}\lambda.
$$
For $t=1$, let  $\sigma_1:\A(S^2_{q1})\to\C$ be the 
$*$-algebra morphism given by the `classical point' \eqref{cpointt1},
$$
\sigma_1(x_0)=1,\qquad\sigma_1(x_1)=0.
$$

\begin{prop}
The top residue of the zeta-type function
$\zeta_a(z):=\Tr_{\HH}(a|D|^{-z})$, with  $a\in\A(S^2_{qt})$, is given by
\begin{equation}                                               \label{eq:rem}
\mathrm{Res}_{z=2}\zeta_a (z)=-\frac{2\mathrm{i}}{\pi}
\int_{S^1}\sigma_t(a)\,\frac{\de\lambda}{\lambda} ,
\qquad  0\leq t\leq 1.
\end{equation}
For $t=1$, $\sigma_1(a)\in\C$ and Equation  \eqref{eq:rem}
simplifies to $\,4\,\sigma_1(a)$.
\end{prop}

\begin{proof}
It is sufficient to prove (\ref{eq:rem}) for the basis elements
$(x_0-t)^jx_1^k$, $j\in\N_0$ and $k\in\Z$, and then extend it to $\A(S^2_{qt})$
by linearity. We use again the notation $x_1^k:=(x_1^*)^{|k|}$ if $k<0$.

Since $\sigma_t\bigl((x_0-t)^jx_1^k\bigr)\propto\delta_{j0}(1-t)^{k/2}\lambda^k$, 
the right hand side of Equation (\ref{eq:rem}) is zero unless \mbox{$j=k=0$}. We next 
show  that the left hand side of (\ref{eq:rem}) also vanishes unless $j=k=0$.
When $j=k=0$, the relation $\zeta_{\id}(z)=4\zeta(z-1)$ fixes the normalization~constant. 

Now
$(x_0-t)^jx_1^k$ is off-diagonal if $k\neq 0$ since it shifts the index $m$ by $k$.
It remains to prove that $\zeta_{(x_0-t)^j}(z)=
\zeta_{\tilde{\varphi}((x_0-t)^j)}(z)\,+$
{\it holomorphic function} has no singularity in $z=2$.
For $j\neq 0$, $\tilde{\varphi}\bigl((x_0-t)^j\bigr)$ satisfies the inequality
$$
\big|\big(v^l_{m,\pm},\tilde{\varphi}\bigl((x_0-t)^j\bigr)v^l_{m,\pm}\big)\big|
\leq c_jq^{l+m},
$$
for some positive constants $c_j$'s. 
From this inequality, we deduce that $\zeta_{\tilde{\varphi}((x_0-t)^j)}(z)$
is a convergent series for all $z$ with $\mathrm{Re}\,z>1$. In particular,
it is finite for $z=2$.
\end{proof}

A natural application of the analysis in this section concerns the construction
of a local representative of the cocycle $\mathrm{ch}^F$, which appeared in
the proof of Proposition \ref{cor:one}.
For our case, Connes-Moscovici theorem~\cite[Theorem~II.3]{CM}  states
that $\mathrm{ch}^F$ is cohomologous to the periodic cyclic cocycle with two components
$(\phi_0,\phi_2)$, given by,
\begin{align*} 
\phi_0(a_0) & =\mathrm{Res}_{z=0}z^{-1}\Tr(\gamma a_0|D|^{-2z})\;\;, \\
\phi_2(a_0,a_1,a_2) & =\mathrm{Res}_{z=0}\Tr(\gamma a_0[D,a_1][D,a_2]|D|^{-2(z+1)})\;\;.
\end{align*}
Smoothing operators do not contribute to $\phi_2$.
Since $[D,a]=\delta\tilde{\varphi}(a)F+\op$, with the map $\tilde{\varphi}$ given in \eqref{tildephi},
we can rewrite $\phi_2$ as,
\begin{align*}
\phi_2(a_0,a_1,a_2) &=\mathrm{Res}_{z=0}\,\Tr\big(\gamma\tilde{\varphi}(a_0)\,\delta
\tilde{\varphi}(a_1)\,\delta\tilde{\varphi}(a_2)|D|^{-2(z+1)}\big)\;.
\end{align*}
Now $\gamma$ is traceless and the remaining operators are diagonal in the
spin index ``$\pm$'';  thus $\phi_2$ is identically zero.
Moreover, $\Tr(\gamma a_0|D|^{-2z})$ is holomorphic for $\mathrm{Re}\,z>1$
where it coincides with the function $\psi(z):=\frac{1}{2}\Tr(\gamma F[F,a_0]|D|^{-2z})$,
which is holomorphic on all $\C$, being $[F,a_0]\in\op$.
Hence,
\begin{equation*}
\phi_0(a_0)=\mathrm{Res}_{z=0}z^{-1}\psi(z)=\psi(0)=\mathrm{ch}^F_0(a_0)\;,
\end{equation*}
with $\mathrm{ch}^F_0$ given by equation (\ref{eq:chzero}).

This shows that, due to the low summability of the Fredholm module, only
the lowest (non-local) component of the periodic cyclic cochain of Theorem~II.3 of \cite{CM} is different from zero, and equals $\mathrm{ch}^F_0$. This fact
was already shown in \cite{DD} for the equatorial Podle\'s sphere.

\section{Final remarks} \label{sec:conclusions}

A crucial ingredient of our analysis was the approximation of 
$\pi(x_i)=\bigl( \pi_{-\half}\oplus\pi_{\half}\bigr) (x_i)$
by the operators $z_i$ defined in \eqref{eq:appr}. 
The proof of Lemma~\ref{appr} shows that a similar approximation
holds true for the representation
$\pi_{-N}\oplus\pi_{N}$ on $\M_{-N}\oplus\M_{N}$, for any $N\in \half \N$.
A careful inspection of the subsequent proofs shows 
that all results are still valid 
and yield corresponding spectral triples with real structure $J$ satisfying
$J^2 = (-1)^N$, and the  pairing of the associated Fredholm module 
with the same projection as 
in Corollary~\ref{cor:one} giving the value $2N$. 
In the classical case $q=1$, these generalized examples
correspond to the quasi-spectral triples studied in \cite{SitarzQuasi}.

\subsection*{Acknowledgements}

It is a pleasure to thank J. Varilly for helpful comments.
F.D'A. L.D. and G.L. thank ESI in Vienna for hospitality.
L.D. was partially supported by the European Union Host
Fellowship for the Transfer of Knowledge Program `Noncommutative
Geometry and Quantum Groups' at Warsaw University.
E.W. was supported by the DFG-grant Wa 1698/2-1.

\end{document}